\documentclass{amsart}

\usepackage{latexsym}
\usepackage{amssymb,amsmath,amsfonts}

\newtheorem{theorem}{Theorem}[section]

\theoremstyle{definition}

\theoremstyle{remark}
\newtheorem{rmk}[theorem]{Remark}

\numberwithin{equation}{section}

\numberwithin{equation}{section}

\begin{document}

\title[Erratum]
 {A note on: Convergence theorems for uniformly L-Lipschitzian Asymptotically Nonexpansive mappings, \emph{Acta Universitatis Apulensis,} no. 20/2009, 183-192.}
 \address{$^1$ Department of Mathematics, Nnamdi Azikiwe University, P. M. B. 5025, Awka, Anambra State, Nigeria.}

 \email{ euofoedu@yahoo.com}
\author[E. U. Ofoedu]{ Eric U. Ofoedu$^{*\;1}$}

\keywords{Nonexpansive mappings, uniformly smooth Banach spaces, viscocity approximation.\\{\indent 2000 {\it Mathematics Subject Classification}. 47H06,
47H09, 47J05, 47J25.\\$^*$The author undertook this work when he visited
 the Abdus Salam International Centre
for Theoretical Physics (ICTP), Trieste, Italy as a visiting fellow.}}

\begin{abstract}
It is our aim in this article to correct the wrong impression created in the paper of A. Rafiq \cite{Rafiq} titled: ``Convergence theorems for uniformly L-Lipschitzian Asymptotically Nonexpansive mappings" which appeared in \emph{Acta Universitatis Apulensis,} no. 20/2009, 183-192.
\end{abstract}

\maketitle
\section{Observations and clarifications.} 

\noindent On page 186 of \cite{Rafiq}, A. Rafiq made the following remark:\\

\noindent{\bf Remark.}\emph{(Remark 3 of \cite{Rafiq}) ``One can easily see that if we take in theorem 3 and 4 $\alpha_n=\frac{1}{n^\sigma}$, $0<\sigma<1,$ then $\displaystyle\sum_{n=1}^\infty \alpha_n=\infty,$ but $\displaystyle\sum_{n=1}^\infty \alpha_n^2=+\infty\;...$" }\\

\noindent This claim false. To see this, consider $\sigma=\frac{2}{3}$. It is easy to see that $0<\frac{2}{3}<1.$ Of course, $\displaystyle \sum_{n=1}^\infty n^{-{\frac{2}{3}}}=+\infty$, but $$\displaystyle \sum_{n=1}^\infty\Big( n^{-{\frac{2}{3}}}\Bigr)^2=\sum_{n=1}^\infty n^{-{\frac{4}{3}}}= \sum_{n=1}^\infty n^{-{(1+\frac{1}{3})}}<+\infty,\;{\rm since}\; p=1+\frac{1}{3}>1.$$ 

\noindent Moreover, Rafiq \cite{Rafiq} claimed that the conditions $$\displaystyle(i)\;\sum_{n\ge 0}\alpha_n^2<+\infty \;\;{\rm and\;\;}(ii)\;\displaystyle\sum_{n\ge 0}\alpha_n(k_n-1)<\infty$$ imposed by Ofoedu \cite{Ofoedu} were dispensed with.\\

\noindent Now, considering the condition $\displaystyle\sum_{n\ge 0}\alpha_n(k_n-1)<\infty$ imposed by Ofoedu \cite{Ofoedu}, it is easy to see that in the paper of Rafiq \cite{Rafiq}, there is no place he/she proved or even assumed that $\displaystyle \sum_{n=1}^\infty (k_n-1)<\infty.$ It is therefore worrisome how on page 190 of \cite{Rafiq} the author got\\

`` so as $j\to \infty$ we have $$\displaystyle \phi(2\phi^{-1}(a_0))\sum_{n=n_0}^\infty b_n\le \|x_0-p\|^2+\sum_{n=n_0}^j(k_n-1)<\infty, "$$ which should in actual sense read \\

`` so as $j\to \infty$ we have $$\displaystyle \phi(2\phi^{-1}(a_0))\sum_{n=n_0}^\infty b_n\le \|x_0-p\|^2+\sum_{n=n_0}^\infty(k_n-1)<\infty ".$$ The rest of the proof of Theorem 8 which is the main theorem of \cite{Rafiq} rested on this.\\

\noindent Even if the condition $\displaystyle \sum_{n\ge 0}(k_n-1)<\infty$ is assumed in \cite{Rafiq}, it is \emph{stronger} than the condition 
$\displaystyle \sum_{n\ge 0}\alpha_n(k_n-1)<\infty$ assumed in \cite{Ofoedu} in the sense that it reduces the class of operators covered in \cite{Ofoedu}. To see this, consider the two sets $$\displaystyle X_1=\Big\{T:K\to K/\;T\;{\rm is\;asympt.\;pseudo. \;with}\;\{k_n\}_{n\ge 1} {\rm \;s.t.\;}\sum_{n\ge 0}(k_n-1)<\infty\Big\}$$
and $$\displaystyle X_2=\Big\{T:K\to K/\;T\;{\rm is\;asympt.\;pseudo. \;with}\;\{k_n\}_{n\ge 1} {\rm \;s.t.\;}\sum_{n\ge 0}\alpha_n(k_n-1)<\infty\Big\}.$$ It is easy to see that $X_1$ \emph{is a proper subset} of  $X_2$. This  follows from the fact that since $\{\alpha_n\}\in (0,1)$ (as assumed in \cite{Ofoedu}), we obtain from \emph{series comparison test} that $\displaystyle \sum_{n=1}^\infty \alpha_n(k_n-1)<\infty$ whenever $\displaystyle \sum_{n=1}^\infty (k_n-1)<\infty$; thus $X_1\subset X_2.$ The converse inclusion is, however, not possible since if $T_1:K\to K$ is an asymptotically pseudocontractive mapping with a sequence $\{k_{1n}\}_{n\ge 1}\subset [1,\infty)$ such that $k_{1n}-1=\frac{1}{\sqrt{n}}\;\forall\;n\in \mathbb{N},$ then $\displaystyle \sum_{n=1}^\infty(k_{1n}-1)=\sum_{n=1}^\infty \frac{1}{\sqrt{n}}=+\infty.$ So, $T_1\notin X_1.$ On the other hand, observe that if we take $\alpha_n=\frac{1}{n}\;\forall\;n\in \mathbb{N},$ then the sequence $\{\alpha_n\}_{n\ge 1}=\Big\{\frac{1}{n}\Big\}_{n\ge 1}$ satisfies the conditions imposed in \cite{Ofoedu}. Furthermore, $$\displaystyle \sum_{n=1}^\infty\alpha_n(k_{1n}-1)=\sum_{n=1}^\infty\Big(\frac{1}{n}.\frac{1}{\sqrt{n}}\Big)=\sum_{n=1}^\infty n^{-{(1+\frac{1}{2}})}<+\infty \;{\rm since\;}1+\frac{1}{2}>1.$$ Thus, $T_1\in X_2$ but $T_1\notin X_1.$ Hence, the class of mappings covered in \cite{Ofoedu} is larger.\\

\noindent  Though the condition the conditions $\displaystyle\sum_{n=1}^\infty b_n=\infty$ and $c_n=o(b_n)$ as imposed on the \emph{iterative parameters} in Theorem 8 of \cite{Rafiq} is milder than the conditions $\displaystyle \sum_{n=1}^\infty \alpha_n=+\infty $  and $\displaystyle\sum_{n\ge 0}\alpha_n^2<+\infty$ imposed in \cite{Ofoedu} (in the sense that the conditions $\displaystyle \sum_{n=1}^\infty \alpha_n=+\infty $  and $\displaystyle\sum_{n\ge 0}\alpha_n^2<+\infty$ implies $\displaystyle\sum_{n=1}^\infty b_n=\infty$ and $c_n=o(b_n)$ with $b_n=\alpha_n$ and $c_n=\alpha_n^2\;\forall\;n\in\mathbb{N}$), it is of interest to note that what is important in construction of iteration processes for approximation of fixed points of nonlinear mappings is to construct iteration processes which are efficient, as simple as possible and in addition accomodates larger classes of nonlinear mappings.\\

\begin{rmk}
Finally, the claims of A. Rafiq \cite{Rafiq} are really strange because none of them made good impact to learning. As shown in \cite{Ofoedu} (see Theorem 3.3 of \cite{Ofoedu}), it is important to note that consideration of the scheme $$x_0\in K,\;x_{n+1}=a_nx_n+b_nT^nx_n+c_nu_n\;n\ge 0$$ (where $\{a_n\}$, $\{b_n\}$, $\{c_n\}$ are sequences in $(0,1)$ with $a_n+b_n+c_n=1$ and $\{u_n\}$ is a bounded sequence in $K$) leads to no further generalization of Theorem 3.2 of \cite{Ofoedu}. The author of \cite{Rafiq} merely repeated the argument of \cite{Ofoedu} and claimed that the results obtained in \cite{Ofoedu} were significantly improved upon. More importantly, his proof of Theorem 8 of \cite{Rafiq} is wrong because its conclusion completely hinged on $\displaystyle\sum_{n\ge 0}(k_n-1)<\infty$ which was neither proved nor stated as part of hypothesis. In the theorems of \cite{Ofoedu}, it should be noted that the strictly increasing continuous fucntion $\phi:[0,+\infty)\to [0,+\infty)$ with $\phi(0)=0$ is further assumed to be onto (surjective) so that $\phi^{-1}$ makes sense. The essence of this paper is to correct all these anomalies.
\end{rmk}

\end{document}